\documentclass[twoside,12pt]{article}
\usepackage{amsfonts,amsmath,amsthm}
\usepackage{tikz,pgfplots}
\usepackage{forest}

\usepackage[a4paper]{geometry}
\newenvironment{itemizeA}
  {\begin{itemize}
    \setlength{\itemsep}{0pt}
    \setlength{\parskip}{0pt}}
  {\end{itemize}}

\theoremstyle{plain}
\newtheorem{theorem}{Theorem}

\newtheorem{proposition}[theorem]{Proposition}

\theoremstyle{definition}
\newtheorem{definition}[theorem]{Definition}
\newtheorem*{example}{Example}
\newtheorem{conjecture}[theorem]{Conjecture}

\theoremstyle{remark}
\newtheorem*{note}{Remark}

\usepackage{accents}
\newlength{\dhatheight}
\newcommand{\doublehat}[1]{%
    \settoheight{\dhatheight}{\ensuremath{\hat{#1}}}%
    \addtolength{\dhatheight}{-0.35ex}%
    \hat{\vphantom{\rule{1pt}{\dhatheight}}%
      \smash{\hat{#1}}}}
\newcommand{\hq}{\hspace*{-1.7 mm}}
\newcommand{\hp}{\hspace*{-2 mm}}
\newcommand{\hs}{\hspace*{-3 mm}}

\def\ff{{\mathbb F}}
\def\ul{\underline}
\def\ds{\displaystyle}
\def\ba{\begin{array}}
\def\ea{\end{array}}
\def\bt{\begin{tabular}}
\def\et{\end{tabular}}
\newcommand \Ci{{\rm C_{I}}}
\newcommand \Cj{{\rm C_{II}}}
\newcommand \Cu{{\rm C_{0}}}
\newcommand \Cv{{\rm C_{1}}}
\def\qq{{\mathbb Q}}
\def\rr{{\mathbb R}}
\def\nn{{\mathbb N}}
\def\dd{{\mathbb D}}
\def\zz{{\mathbb Z}}

\begin{document}
{\bf\noindent  {\large V Tree} --- \\
Continued Fraction Expansion,
Stern-Brocot Tree,  Minkowski's $?({\bf x})$ Function\\
In Binary: Exponentially  Faster
}  
  \vspace*{ 5mm}

  \noindent
  
\noindent   Michael Vielhaber
  
  \vspace*{ 5mm}
  
\noindent   Hochschule Bremerhaven, {FB2},
  An der Karlsburg 8,
\noindent   {D--$27568$}  Bremerhaven, Deutschland\\
{\tt vielhaber@gmail.com}

  \vspace*{ 12mm}

{\bf Abstract}\quad
  The Stern-Brocot tree and  Minkowki's question mark function $?(x)$ (or Conway's box function) are related to the continued fraction expansion of numbers from $\qq$ with unary encoding of the partial denominators.

  We first define binary encodings $\Ci,\Cj$ of the natural numbers, adapted to the Gau\ss-Kuz'min measure for the distribution of partial denominators.

  Then we define the V$_1$ tree as analogue to the Stern-Brocot tree, using the binary encondings 
  $\Ci,\Cj$.
  We shall see that all numbers with denominator $q$ are present in the first $3.44\log_2(q)$ levels, instead of $1/q$ appearing in level $q$ in the Stern-Brocot tree.
  The extension of the V$_1$  tree, the V tree, covers all numbers from $\qq$ exactly once.
  We also define the binary version of Minkowski's question mark function, $?_V$, and conjecture that it has no derivative at rational points (for the original, $?'(x)=0, x\in\qq_1$).

\vspace*{5mm}

\centerline{
  {\it ``Read the classics''} \ \ ---\ \ Edwards \cite[S.~ix]{Edwa}
}

\vspace*{5mm}

\begin{table}[!b]
  \rule{6 cm}{0.4 pt}
  
  \noindent {\bf Notation:}
  
  \noindent $\nn= \{1,2,3,\dots\}, \nn_0 = \{0,1,2,3,\dots\}$
  
  \noindent $\dd := \{a / 2^k\colon a\in\zz \mbox{\rm\ odd}, k\in\nn_0\}$, dyadic fractions
  
  \noindent For $X\subset \rr$, $X_1 := X\cap (0,1), X^+ := X\cap (0,\infty)$: $\dd_1,\qq_1,\rr_1$,\  $\dd^+,\qq^+,\rr^+$
  
  \noindent $A=\{0,1\}$ is the binary alphabet
  
  \noindent $A^*=\{\varepsilon,0,1,00,01,10,11,000,\dots\}$ and  $A^\omega$ are  the finite, resp.~infinite words over $A$
  
  \noindent For $v= v_1v_2\dots v_{|v|}\in A^*$: $(v)_2 = \sum_{k=0}^{|v|-1} v_{|v|-k} 2^k\in\nn_0$ in binary, $(\varepsilon)_2=0$
  
  \noindent 0/1-inversion: $\overline 0 = 1,\ \overline 1 = 0,\ \overline{10011} = 01100$
  
  \noindent $\varphi_b := (b+\sqrt{b^2+4})/2$ is the larger root of $x^2= x\cdot b+1$, eigenvalue/growth rate of $A_k, B_k$ if all  PDs $b_i=b$ $(\lambda_b$ in \cite{Dush}) \quad $\varphi_1=1.618, \varphi_2=2.414, \varphi_4=4.236$ 
\end{table}

\noindent {\bf Introduction / Motivation:}
In the theory of stream ciphers, the continued fraction expansion of formal power series from $\ff_2[[x^{-1}]]$ leads to an isometry between the coefficient series $s = (s_k)$ and the encoding of the partial denominators (which are polynomials from $\ff_2[x]$ in this case),
${\mathbb K}\colon \ff_2^\omega\ni s\mapsto d\in\ff_2^\omega$ is an isometry, 
  $${\mathbb K}\colon s\mapsto G(s) =\sum_{k\in\nn}s_kx^{-k} = [b_1,b_2,b_3,\dots] \mapsto C_{\ff_2[x]}(b_1)|C_{\ff_2[x]}(b_2)|C_{\ff_2[x]}(b_3)|\dots = (d_k).$$
  
  This  aestetically pleasing result motivated the paper, considering the same problem for  $\rr$.
  Since $\rr$ is Archimedean, while $\ff_2^\omega, \ff_2[[x^{-1}]]$ are ultrametric, we shall see (Gau\ss-Kuz'min measure) that an exact isometric result is impossible (with the exception of an ideal context-sensitive ``L\'evy encoding'', which would though just be the identity on $A^\omega$).

  Nevertheless, the presented encodings are a huge step forward in terms of the expected codeword length $H$:
  From  $H_{SB}=\infty$ for the unary, Stern-Brocot case, to $H_{C_I,C_{II}} = 3.507$ for our codes, near the optimum  $H_{L\acute{e}vy}=3.423$.
  For more implementation details see~\cite{Vseta2020}.
  
   See Berlekamp \cite{Berl} and Massey \cite{Mass} for the general solution, Dornstetter \cite{Dorn} and V.~\cite{Vieee} for the isometric adaptation, Niederreiter and V.~\cite{NV}, Canales and V.~\cite{Cana} for applications.

\subsection*{I -- Definitions}

\begin{definition} Binary Encodings $\Ci, \Cj$ (see Table~\ref{fig:Codes})
  \label{binary}
  
  Let $\Ci\colon \nn\cup \{\aleph_0\} \to A^*\cup \{0^\omega\}$
  $$b = \sum_{k=0}^l b_k2^k\ \mapsto \ \Ci(b) = 0^{l}1\overline b_{l-1}\overline b_{l-2}\dots \overline b_1\overline b_0,\hspace*{17mm}$$
  where  $l = \lfloor \log_2(b)\rfloor$, and $\Ci(\aleph_0) :=  0^\omega$,  be a complete prefixfree code.

  Let $\Cj\colon \nn\cup \{\aleph_0\} \to A^*\cup \{1^\omega\}$
  $$b = \sum_{k=0}^l b_k2^k\  \mapsto \  \Cj(b) = 1^{l}0b_{l-1}b_{l-2}\dots b_1b_0 = \overline{\Ci(b)},$$
  with $l$ as before and $\Cj(\aleph_0) :=  1^\omega$,  be the 0/1-inverse of $\Ci$, again a complete prefix code.

We have $l_{I,II}(b): = 1 + 2\cdot \log_2(b)$ as length of the codewords for $b$.
\end{definition}

\begin{table}[!h]
  \begin{center}
    {\small{
        \begin{tabular}{rllccc}
          $b$&$\Ci(b)$&$\Cj(b)$&$l_{\rm I,II}(b)$&$l_{GK}$&$\mu_{GK}$\\
          \cline{1-6}
          1 &1        &0         &1& 1.269& 0.4150\\
          \cline{1-4}
          2 &011      &100       &3& 2.557& 0.1699\\
          3 &010      &101       &3& 3.425& 0.0931\\
          \cline{1-4}                
          4 &00111    &11000     &5& 4.086& 0.0588\\
          5 &00110    &11001     &5& 4.621& 0.0406\\
          6 &00101    &11010     &5& 5.071& 0.0297\\
          7 &00100    &11011     &5& 5.460& 0.0227\\
          \cline{1-4}                
          8 &0001111  &1110000   &7& 5.802& 0.0179\\ 
          9 &0001110  &1110001   &7& 6.108& 0.0144\\
          10&0001101  &1110010   &7& 6.384& 0.0119\\
          11&0001100  &1110011   &7& 6.636& 0.0100\\
          12&0001011  &1110100   &7& 6.868& 0.0085\\
          13&0001010  &1110101   &7& 7.082& 0.0073\\
          14&0001001  &1110110   &7& 7.282& 0.0064\\
          15&0001000  &1110111   &7& 7.468& 0.0056\\
          \cline{1-4}
          16&000011111&111100000 &9& 7.644& 0.0050\\  
          \vdots&&&&\\ 
          31&000010000    &111101111     & 9& 9.471& 0.0014\\
          \cline{1-4}                              
          32&00000111111  &11111000000   &11& 9.559 & 0.0013\\
          \vdots&&&&\\ 
          63&00000100000  &11111011111   &11& 11.471& 0.00035\\
          \cline{1-4}                
          64&0000001111111&1111110000000 &13& 11.516& 0.00034\\
          \vdots&&&\\  
          $\aleph_0$&$0^\omega$&$1^\omega$&---&---&---\\
        \end{tabular}              
    }}
    
    \vspace*{3 mm}
    
    (the columns $l_{GK}$ and $\mu_{GK}$ are explained in Theorem~\ref{GK})
  \end{center}

  \caption{\label{fig:Codes} Codes $\Ci$ and $\Cj$  for partial denominators.}
\end{table}

\begin{example}
  $b = 14 = (1\framebox{110})_2 \mapsto 000\ul 1\framebox{001} = \Ci(b)$ and $111\ul 0\framebox{110} = \Cj(b)$
\end{example}

\begin{definition}  Binary V question mark functions  $?_V$ and $?_V^{-1}$
  \label{def:QMF}

  $(i)$ For {\large{$\frac{p}{q}$}} $\in\qq_1$, let {\large{$\frac{p}{q}$}} $= [b_1,b_2,\dots,b_{2l}]$ be its continued fraction expansion (CFE) with an even number of partial denominators (PD)  (see Appendix 1).
  
  We define the function $?_V\colon \qq_1\to A^*$ by 
  $$?_V\left(\frac{p}{q}\right) := \Ci(b_1)|\Cj(b_2)|\dots|\Ci(b_{2l-1})|\Cj(b_{2l})\quad  \backslash\  10^*.$$
  The operation $\backslash 10^*$ removes all, if any, trailing zeroes and then one symbol 1.
  This affects at most $\Cj(b_{2l})$ and, only in case of $b_{2l}=1, \Cj(1)=0$, also affects $\Ci(b_{2l-1})$.

  $(ii)$ For any $v\in A^*$, extended to the infinite word $v10^\omega$, let 
  $$v10^\omega =  \Ci(b_1)|\Cj(b_2)|\dots|\Ci(b_{2l-1})|\Cj(b_{2l})|\Ci(\aleph_0)$$
  be the decomposition of $v10^\omega$ into encodings, starting with $\Ci$.
  Since $\Ci, \Cj$ are complete and prefixfree, this is always possible, in a unique way.

  Then $?_V^{-1}\colon A^* \to \qq_1$ is defined by 
  $$?_V^{-1}(v) := \frac{p}{q} = [b_1,b_2,\dots,b_{2l}].$$

  By construction, we have $?_V^{-1}(?_V(p/q)) = p/q$  and  $?_V(?_V^{-1}(v)) = v$  for all $v\in A^*$ and $p/q\in \qq_1$.

  $(iii)$ We define real-valued functions $\overline ?_V$ and $\overline ?_V^{-1}$ from $\rr_1$ to $\rr_1$ by first defining
  $$\overline ?_V\colon\qq_1\to \dd_1,\quad \overline ?_V(p/q) := \iota_{AD}\left(?_V(p/q)\right)$$
  with $\iota_{AD}$ from Appendix~2 on the equivalence $\nn\equiv A^* \equiv \dd_1$, and then
  $\overline ?_V\colon \rr_1\to \rr_1$ by continuous extension. Also, first   
  $$\overline ?^{-1}_V\colon\dd_1\to \qq_1,\quad \overline ?^{-1}_V(d) := ?^{-1}_V(\iota_{DA}(d))$$
  and then continuously extending to $\rr_1$.
\end{definition}

\begin{definition} V$_{10}$ tree for $\qq_1$
    
  We define the  V$_{10}$ tree as an infinite binary tree with label $?_V^{-1}(v)\in\qq_1$ at the node with symbolic address $v$ (see Appendix~4 on trees and addresses).

\end{definition}

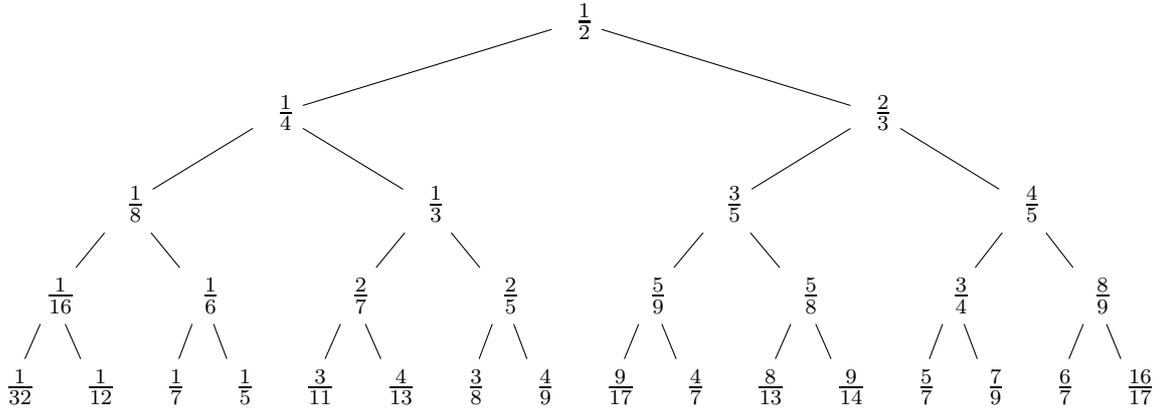
\begin{figure}[!t]
  \begin{center}
    {\small{
        \begin{forest}
          [$\ba{c}\frac{1}{2}\ea\hp$
            [$\ba{c}\frac{1}{4}\ea\hp$
              [$\ba{c}\frac{1}{8}\ea\hp$
                [$\ba{c}\frac{1}{16}\ea\hp$
                  [$\ba{c}\frac{1}{32}\ea\hp$]
                  [$\ba{c}\frac{1}{12}\ea\hp$]
                ]
                [$\ba{c}\frac{1}{6}\ea\hp$
                  [$\ba{c}\frac{1}{7}\ea\hp$]
                  [$\ba{c}\frac{1}{5}\ea\hp$]
                ]
              ]
              [$\ba{c}\frac{1}{3}\ea\hp$
                [$\ba{c}\frac{2}{7}\ea\hp$
                  [$\ba{c}\frac{3}{11}\ea\hp$]
                  [$\ba{c}\frac{4}{13}\ea\hp$]
                ]
                [$\ba{c}\frac{2}{5}\ea\hp$
                  [$\ba{c}\frac{3}{8}\ea\hp$]
                  [$\ba{c}\frac{4}{9}\ea\hp$]
                ]
              ]
            ]
            [$\ba{c}\frac{2}{3}\ea\hp$
              [$\ba{c}\frac{3}{5}\ea\hp$
                [$\ba{c}\frac{5}{9}\ea\hp$
                  [$\ba{c}\frac{9}{17}\ea\hp$]
                  [$\ba{c}\frac{4}{7}\ea\hp$]
                ]
                [$\ba{c}\frac{5}{8}\ea\hp$
                  [$\ba{c}\frac{8}{13}\ea\hp$]
                  [$\ba{c}\frac{9}{14}\ea\hp$]
                ]
              ]
              [$\ba{c}\frac{4}{5}\ea\hp$
                [$\ba{c}\frac{3}{4}\ea\hp$
                  [$\ba{c}\frac{5}{7}\ea\hp$]
                  [$\ba{c}\frac{7}{9}\ea\hp$]
                ]
                [$\ba{c}\frac{8}{9}\ea\hp$
                  [$\ba{c}\frac{6}{7}\ea\hp$]
                  [$\ba{c}\frac{16}{17}\ea\hp$]
                ]
              ]
            ]
          ]
        \end{forest}
    }}                             
  \end{center}                     
  \caption{\label{fig:VTpos2} V$_{10}$ tree on $\qq_1$.}
\end{figure}

\begin{definition}
  \label{def:QMFhat}
  V$_{1}$ tree and V question mark functions $\hat{?}_V,\hat{?}^{-1}_V$ for $\qq^+$ 

  From $\qq_1$ to $\qq^+$ by multiplicative inversion:

  $(i)$ For $m\in A, v\in A^*$, let
  \begin{eqnarray*}
  \hat{?}^{-1}\colon A\times A^*&\to& \qq^+\\
  \hat{?}^{-1}(0v)&=&{?}^{-1}(v)\\
  \hat{?}^{-1}(1v)&=&\left({?}^{-1}(\overline v)\right)^{-1}
  \end{eqnarray*}
  where $\overline v$ is the 0/1-inverted address.

  We also set $\hat{?}^{-1}(\varepsilon) = 1$.
  Then $\hat{?}^{-1}\colon A^*\to \qq^+$ is defined on all of $A^*$.

  $(ii)$ We define the V$_1$ tree as an infinite binary tree with label $\hat{?}^{-1}(v)$ at the node with address~$v$.

  $(iii)$ Let $\hat{?}\colon \qq^+\to A^*$ be the inverse function to $\hat{?}^{-1}$,
  $$\hat{?}\left(\frac{p}{q}\right) =\left\{
  \ba{ll}
  0|?_V(\frac{p}{q}),& p<q,\\
\varepsilon,&p=q,\quad {\it i.e.\ }p/q=1,\\
  1|\overline{?_V(\frac{q}{p})},& p > q.
  \ea\right.$$
  
\end{definition}

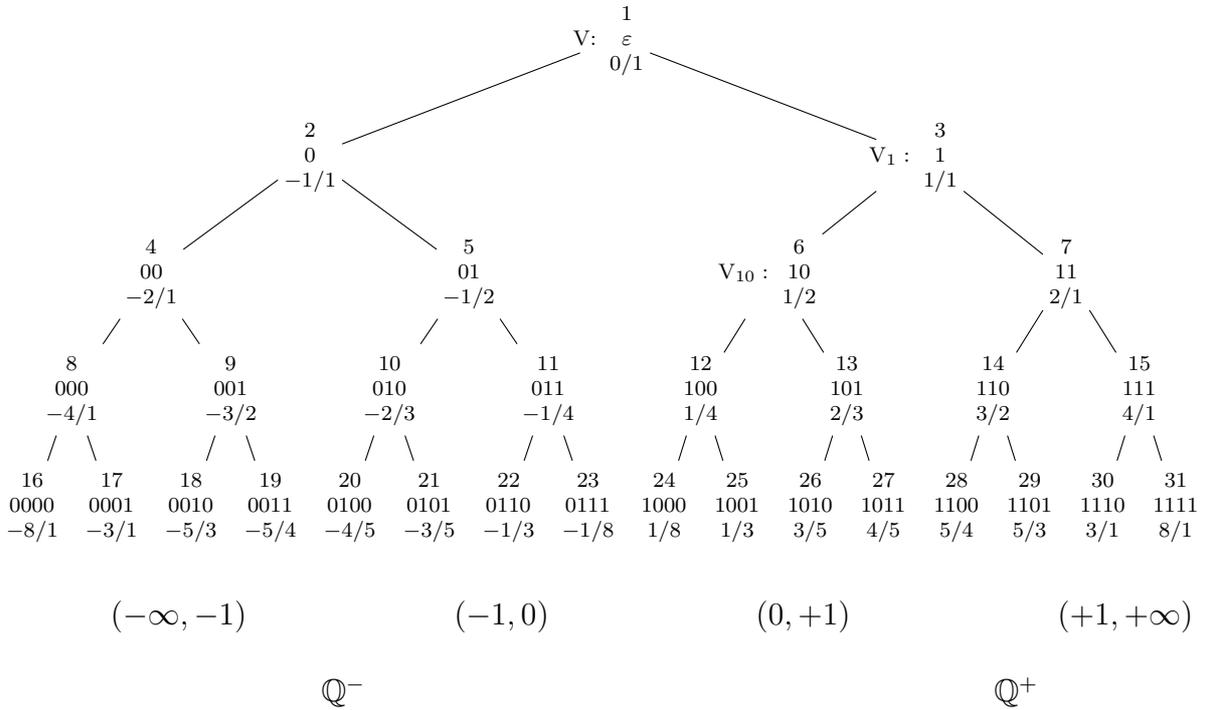
\begin{figure}[!h]
  \begin{center}
    {\scriptsize{
        \begin{forest}
          [{\hp V:$\ba{c}1\\\varepsilon\\{0}/{1}\ea$\hp}
            [{\hp$\ba{c}2\\0\\{-1}/{1}\ea$\hp}
              [{\hp$\ba{c}4\\00 \\{-2}/{1}\ea$\hp}
                [{\hp$\ba{c}8\\000 \\{-4}/{1}\ea$\hp}
                  [{\hp$\ba{c}16\\0000 \\{-8}/{1}\ea$\hp}]
                  [{\hp$\ba{c}17\\0001 \\{-3}/{1}\ea$\hp}]
                ]
                [{\hp$\ba{c}9\\001 \\{-3}/{2}\ea$\hp}
                  [{\hp$\ba{c}18\\0010 \\{-5}/{3}\ea$\hp}]
                  [{\hp$\ba{c}19\\0011 \\{-5}/{4}\ea$\hp}]
                ]
              ]
              [{\hp$\ba{c}5\\01 \\{-1}/{2}\ea$\hp}
                [{\hp$\ba{c}10\\010 \\{-2}/{3}\ea$\hp}
                  [{\hp$\ba{c}20\\0100 \\{-4}/{5}\ea$\hp}]
                  [{\hp$\ba{c}21\\0101 \\{-3}/{5}\ea$\hp}]
                ]
                [{\hp$\ba{c}11\\011 \\{-1}/{4}\ea$\hp}
                  [{\hp$\ba{c}22\\0110 \\{-1}/{3}\ea$\hp}]
                  [{\hp$\ba{c}23\\0111\\{-1}/{8}\ea$\hp}]
                ]
              ]
            ]
            [{\hp V$_1:$$\ba{c}3\\1\\{1}/{1}\ea$\hp}
              [{\hp V$_{10}:$$\ba{c}6\\10 \\{1}/{2}\ea$\hp}
                [{\hp$\ba{c}12\\100 \\{1}/{4}\ea$\hp}
                  [{\hp$\ba{c}24\\1000 \\{1}/{8}\ea$\hp}]
                  [{\hp$\ba{c}25\\1001 \\{1}/{3}\ea$\hp}]
                ]
                [{\hp$\ba{c}13\\101 \\{2}/{3}\ea$\hp}
                  [{\hp$\ba{c}26\\1010 \\{3}/{5}\ea$\hp}]
                  [{\hp$\ba{c}27\\1011 \\{4}/{5}\ea$\hp}]
                ]
              ]
              [{\hp$\ba{c}7\\11\\{2}/{1}\ea$\hp}
                [{\hp$\ba{c}14\\110 \\{3}/{2}\ea$\hp}
                  [{\hp$\ba{c}28\\1100\\{5}/{4}\ea$\hp}]
                  [{\hp$\ba{c}29\\1101 \\{5}/{3}\ea$\hp}]
                ]
                [{\hp$\ba{c}15\\111 \\{4}/{1}\ea$\hp}
                  [{\hp$\ba{c}30\\1110 \\{3}/{1}\ea$\hp}]
                  [{\hp$\ba{c}31\\1111\\{8}/{1}\ea$\hp}]
                ]
              ]
            ]
          ]
        \end{forest}
    }}                             
  \end{center}                     
  
  \hspace*{1 mm}
  \hfill
  $(-\infty,-1)$
  \hfill
  \hfill
  $(-1,0)$
  \hfill
  \hfill
  $(0,+1)$
  \hfill
  \hfill
  $(+1,+\infty)$
  \hspace*{1 mm}
  \\
  \\
  \hspace*{1 mm}
  \hfill
  \hfill
  $\qq^-$
  \hfill
  \hfill
  \hfill
  \hfill
  $\qq^+$
  \hfill
  \hspace*{1 mm}
  
  \caption{\label{fig:VTree} Full V Tree on $\qq$ with  addresses.}
\end{figure}

\begin{definition} V tree and V question mark functions $\doublehat{?}_V,\doublehat{?}^{-1}_V$ for $\qq$ 
  \label{def:QMFdhat}

  From $\qq^+$ to $\qq$ by additive inversion:

  $(i)$ For $a\in A, v\in A^*$, let
  \begin{eqnarray*}
    \doublehat{?}^{-1}\colon A\times A^*&\to& \qq\\
    \doublehat{?}^{-1}(0v)&=&-\left(\hat{?}^{-1}(\overline v)\right),\\
    \doublehat{?}^{-1}(1v)&=&\hat{?}^{-1}(v).\\
  \end{eqnarray*}
  We also set $\doublehat{?}^{-1}(\varepsilon) = 0$.
  Then $\doublehat{?}^{-1}\colon A^*\to \qq$ is defined on all of $A^*$.

  $(ii)$ We now  define  the V tree (see Figure~\ref{fig:VTree}) as an infinite binary tree with label $\doublehat{?}^{-1}(v)$ at the node with address $v$.

  $(iii)$ Let $\doublehat{?}\colon \qq\to A^*$ be the inverse function to $\doublehat{?}^{-1}$,

  $$\doublehat{?}\left(\frac{p}{q}\right) =\left\{
  \ba{ll}
  0|\overline{\hat{?}_V(-\frac{p}{q})},& p/q < 0,\\
  \varepsilon,&p/q = 0,\\
  1|\hat{?}_V(\frac{p}{q}),& p/q > 0.\\
  \ea\right.$$
  
\end{definition}

\begin{definition} Sequences V, V$_1$, V$_{10}$
  
  Reading out the values from the V, V$_1$, and V$_{10}$ trees in the order of the numerical addresses (breadth first), we obtain the following 3 sequences:

  \vspace*{ 1mm}
  
  \noindent\hspace*{2mm}$ V = $ {\large{$ \left(\frac{p_n}{q_n}\right)_{n=1}^\infty =
      \left(\frac{0}{1},\frac{-1}{1},\frac{1}{1},
      \frac{-2}{1},\frac{-1}{2},\frac{1}{2},\frac{2}{1},
      \frac{-4}{1},\frac{-3}{2},\frac{-2}{3},\frac{-1}{4},\frac{1}{4},\frac{2}{3},\frac{3}{2},\frac{4}{1},
      \frac{-1}{8},
      \dots\right)$ }}\ $ \equiv \qq$

  \vspace*{ 1mm}
  
  \noindent\hspace*{2mm}$ V_1 =$ {\large{$\left(\frac{1}{1},
      \frac{1}{2},\frac{2}{1},
      \frac{1}{4},\frac{2}{3},\frac{3}{2},\frac{4}{1},
      \frac{1}{8},\frac{1}{3},\frac{3}{5},\frac{4}{5},
      \frac{5}{4},\frac{5}{3},\frac{3}{1},\frac{8}{1},
      \frac{1}{16},\frac{1}{6},\frac{2}{7},\frac{2}{5},
      \frac{5}{9},\frac{5}{8},\frac{3}{4},\frac{8}{9},
      \frac{9}{8},\frac{4}{3},
      \dots\right)$}} $\  \equiv \qq^+$

  \vspace*{ 1mm}
  
  \noindent\hspace*{2mm}$ V_{10} =$   {\large{$
      \left(\frac{1}{2},\frac{1}{4},\frac{2}{3},\frac{1}{8},\frac{1}{3},\frac{3}{5},\frac{4}{5},\frac{1}{16},\frac{1}{6},
      \frac{2}{7},\frac{2}{5},
      \frac{5}{9},\frac{5}{8},\frac{3}{4},\frac{8}{9},\frac{1}{32},\frac{1}{12},\frac{1}{7},\frac{1}{5},
      \frac{3}{11},\frac{4}{13},\frac{3}{8},\frac{4}{9},
      \dots\right)$}}$\ \ \  \equiv \qq_1.$

  \vspace*{ 1mm}
  
  \noindent We show in Theorem~\ref{VisQ} that indeed the sequences V, V$_1$, V$_{10}$ are a  complete ordering of all elements of  $\qq,\qq^+, $ and $\qq_1$, respectively, each element appearing exactly once.
\end{definition}

\subsection*{II -- Rationale}

\begin{definition} Unary Encodings
  \label{unary}
  
  $(i)$ Let
  $$\Cu\colon \nn\cup \{\aleph_0\} \to A^*\cup \{0^\omega\},\quad \Cu(b) = 0^{b-1}1,\quad \Cu(\aleph_0) =  0^\omega,$$
  $$\Cv\colon \nn\cup \{\aleph_0\} \to A^*\cup \{1^\omega\},\quad \Cv(b) = 1^{b-1}0,\quad \Cv(\aleph_0) =  1^\omega$$
  be complete prefixfree codes.

  $(ii)$ Let
  $$\Cu'\colon \nn\cup \{\aleph_0\} \to A^*\cup \{0^\omega\},\quad \Cu'(b) = 0^{b},\quad \Cu'(\aleph_0) =  0^\omega,$$
  $$\Cv'\colon \nn\cup \{\aleph_0\} \to A^*\cup \{1^\omega\},\quad \Cv'(b) = 1^{b},\quad \Cv'(\aleph_0) =  1^\omega$$
  be codes, which are, however,  neither complete nor  prefixfree.
\end{definition}

\begin{theorem} $($Gau\ss{}-Kuz'min-Khinchin-L\'evy {\rm\cite{Khin}\cite{Kuzm}\cite{Levy})}\label{GK}
  
  For almost all  values $r\in\rr$, we have:
  
  $(i)$ The probability for a partial denominator $b\in\nn$, its Gau\ss-Kuz'min measure,  is
  $$\mu_{GK}(b) = -\log_2\left(1-\frac{1}{(b+1)^2}\right) = \log_2\left(1+\frac{1}{b(b+2)}\right).$$
  
  $(ii)$ The geometric average of the partial denominators is  Khinchin's constant
  $$K := \lim_{n\to\infty}\sqrt[n]{b_1\cdot b_2\cdots b_n} = 2.68545.$$
  
  $(iii)$ The average gain in precision, per partial denominator in bits, is\\
  \centerline{$\displaystyle \frac{\pi^2}{6\ln(2)^2} = 3.42371\dots{} =: H_{L\acute{e}vy},$  with $2^{3.42371/2}= 3.27582\dots$ being L\'evy's constant.}
\end{theorem}

One obtains Khinchin's constant as geometric average over the Gau\ss-Kuz'min measure,
$$K = \prod_{b\in\nn} b^{\mu_{GK}(b)}.$$

\begin{note}
  The exceptions to this result are
  \begin{itemizeA}
  \item[---] rational numbers  (Euclid, \cite[Liber VII, \S 1+2]{Eukl})
  \item[---] quadratic-algebraic numbers (Lagrange \cite{Lagr}
  \item[---] powers $e^{2/k}$ (Euler \cite{Eule}, Hurwitz \cite{Hurw})
  \item[---] Liouville numbers (Maillet \cite{Mail}, Liouville \cite{Liou})
  \item[---] numbers with bounded PDs (Shallit \cite{Shal}, Jenkinson \cite{Jenk})
  \end{itemizeA}
  and some more, altogether  a set of positive Hausdorff dimension, but measure zero.
  The encodings $\Ci, \Cj$ are modelled along the Gau\ss-Kuz'min measure, which suggests ``ideal'' codeword lengths $l_{GK} = -\log_2(\mu_{GK}(b))$ (see Table~\ref{fig:Codes}), which are however non-integral. 
\end{note}

\begin{note} {Codeword lengths}\quad
  Besides $\mu_{GK}$ from Theorem~\ref{GK}(i) with (non-integral) codeword length $l_{GK}(b) := -\log_2(\mu_{GK}(b)$, we need $l_{I,II}(b) := 1+2\log_2(b)$ from Definition~\ref{binary},  and $l_{SB}(b) := b$ for the unary encoding from Definition~\ref{unary}.
\end{note}

\begin{proposition}   Average Codelength
  \label{entro}

  For each code $X$ with codeword length $l_X$, we define the average codeword length or entropy under the Gau\ss-Kuz'min distribution as 
  $$H_X := \sum_{b\in\nn} l_X(b)\cdot \mu_{GK}(b).$$
  Also, let $H_{L\acute{e}vy} = 3.42371\dots$ from Theorem~{\rm \ref{GK}(iii)} as lower bound.

  By numerical evaluation, we obtain the average codeword lengths/entropies $($Table~{\rm\ref{fig:Entropy}}$)$.

\begin{table}[!h]
  \begin{center}
    {\small{
        \begin{tabular}{l|cccc}
          &\multicolumn{4}{c}{Average codeword length}\\
          Code $X$ &L\'evy&GK         & ${\Ci,\Cj}$ &{SB}\\
          \cline{1-5}
          $H_X$ & 3.42371       & { 3.43246}& { 3.50698}  & { $\infty$}\\
        \end{tabular}              
    }}                             
  \end{center}                     
  \caption{\label{fig:Entropy} Entropies resulting from various codes.}
\end{table}

\end{proposition}
      
\begin{theorem} The Stern-Brocot tree and Minkowski's question mark function $?(x)$ 

  $(i)$ Let $?, ?^{-1}$ be defined analogously to Definition~{\rm \ref{def:QMF}}, but using codes $\Cu, \Cv$.\\
  Then $?(x)$ is Minkowski's question mark function $($equal to Conway's box function {\rm\cite{Conw})}.
  $$?(p/q) = \iota_{AD}(\Cu(b_1)|\Cv(b_2)|\dots|\Cv(b_{2l}).$$
  
  $(ii)$  Let $\hat{?}, \hat{?}^{-1}$ be defined analogously to Definition~{\rm \ref{def:QMFhat}}, but using codes $\Cu, \Cv$.\\
  The infinite binary tree with $\hat{?}^{-1}(v)$ as label at node $v$ then is the Stern-Brocot tree $($Stern~{\rm\cite{Ster}}, Brocot~{\rm\cite{Broc})}, see Appendix~$4$.
\end{theorem}

\begin{proof}
  $(i)$ For $p/q = [b_1,b_2,\dots,b_{2l}]$, we have

  $$?(p/q) = \iota_{AD}(\Cu(b_1)|\Cv(b_2)|\dots|\Cv(b_{2l})$$
  $$ = 2\cdot \iota_{AD}(\Cu'(b_1)|\Cv'(b_2)|\dots|\Cv'(b_{2l}))  =  2\cdot\sum_{k=1}^l (-1)^{k+1} 2^{-\sum_{i=1}^k b_i}.$$ 

  The last  representation for $?(x)$ was introduced by Denjoy~\cite{Denj}.

 $(ii)$ follows from $(i)$ and the known correspondence between  Minkowski's $?(x)$ function and the entries in the Stern-Brocot tree.\end{proof}

\begin{note}
  The V$_1$  tree is the analogue to the Stern-Brocot tree for binary encoded partial denominators.  
  More on the Stern-Brocot tree and Minkowski's $?(x)$ function can be found in Salem~\cite{Sale}, Viader~{\it et al.}~\cite{Viad}, and Vepstas \cite{Veps}.
\end{note}

\subsection*{III -- Properties}

\begin{theorem} Equivalence of the V$_{10}$, V$_1$, V  sequences with $\qq_1,\qq^+,\qq$
  \label{VisQ}
  
  $(i)$ The sequence V$_{10}$ contains every element from $\qq_1$ exactly once.

  $(ii)$ The sequence V$_{1}$ contains every element from $\qq^+$ exactly once.

  $(iii)$ The sequence V contains every element from $\qq$ exactly once.

\end{theorem}

\begin{proof}

  $(i)$ Every $p/q\in\qq_1$ has a finite CFE $[b_1,\dots,b_{2l}]$ (where $b_0=0$ can be omitted) with encoding
  $C(p/q)= \Ci(b_1)|\dots|\Cj(b_{2l})$ and resulting address $v = C(p/q)\backslash 10^*$.

  Hence, $p/q$ is present in the V$_{10}$ tree and the V$_{10}$ sequence at node $v$ and place $n$, respectively,  with $n=(1v)_2$, and only there.
  Different $p/q$ lead to different CFEs, since   $\Ci, \Cj$ are prefixfree, and the operation $\backslash 10^*$ keeps the node addresses  $v$ different.

  $(ii)$ By construction, the V$_1$ tree and thus the V$_1$ sequence contain exactly once every element from $\qq_1$ (by $(i)$), in the left subtree, their multiplicative inverses in the right subtree, and 1 as root or first element, respectively.
  
  Since for every $p/q\in\qq^+$, we either have\\
  \hspace*{2cm}
  \bt{l}
  $p<q$ and thus $p/q \in\qq_1$, or\\
  $p>q$ and thus $q/p \in\qq_1$, or\\
  $p=q$ and thus $p/q = 1$,\\
  \et

\noindent and these cases are mutually exclusive, we are done.

  $(iii)$ By construction, the V tree and thus the V sequence contain exactly once every element from $\qq^+$ (by $(ii)$), in the right  subtree, their additive inverses in the left subtree, and 0 as root or first element, respectively.
  
  Since for every $p/q\in\qq$, we either have\\
  \hspace*{2cm} 
  \bt{l}
  $p/q>0$ and thus $+p/q \in\qq^+$, or\\
  $p/q<0$ and thus $-p/q \in\qq^+$, or\\
  $p/q=0$,\\
  \et
  
\noindent   and these cases are mutually exclusive, we are done.
\end{proof}

\begin{theorem} Monotonicity of Codes $\Ci, \Cj$ 
  
  Let $r=[b_1,b_2,\dots] ,r'=[b'_1,b'_2,\dots]\in \rr_1$ with $r < r'$, and
  $$C(r)=\Ci(b_1)|\Cj(b_2)|\dots,\quad  C(r')=\Ci(b'_1)|\Cj(b'_2)|\dots{} \in A^\omega$$
  their  encodings $($for $r,r'\in\qq$, terminate with $\Ci(\aleph_0) = 0^\omega)$.

  Then $C(r) < C(r')$ in lexicographic order.
\end{theorem}

\begin{proof}
  Let $C(r)$ and $C(r)$ have identical encodings $\Ci(b_1)=\Ci(b'_1), \Cj(b_2)=\Cj(b'_2), \dots$ until the first $b_l\neq b'_l$.
  
  $(i)$ If $l$ is odd, $b_l > b_l'$ implies
  $r < r'$, regardless of the further PDs, see \cite[Satz 2.9]{Perr}.
  Also, $b_l > b_l'$ implies $\Ci(b_l) < \Ci(b'_l)$ in lexicographical order.
  For $b_l < b'_l$, all relations are inverted.
  We thus obtain that $r < r'$ implies $C(r) < C(r')$.
  
  $(ii)$ If $l$ is even, $b_l > b_l'$ implies
  $r > r'$ and  $\Ci(b_l) > \Ci(b'_l)$.
  Again for $b_l < b'_l$, all relations are inverted.
  Therefore, also for even $l$,  $r < r'$ implies $C(r) < C(r')$.
\end{proof}

\begin{note}
  We are entering the realm of ``Experimental Mathematics''

\nopagebreak
  
  The following two theorems have been proved (or ``proved'') by verifying $2^{30}$ cases (nodes of the respective tree).
  The author sees no chance of changing circumstances in levels 31 and below, in view of Appendix~3.
  (Consult \cite{ExMath} for philosophical consolation :-)
\end{note}

\begin{theorem} $[$Conjecture$]$  Determinants between Neighbour Nodes
  \label{pqpq}
  
  Let the $2^n-1$ values in levels $1,\dots, n$ be linearized, {\it i.e.}~starting with the root, we place the $2^l-1$ elements from levels $1,\dots, l$ between the $2^l$ elements in level $l+1$, for $l = 1,2,\dots, n-1$.

  Let then $p_k/q_k$ be the value in position $k, 1\leq k\leq 2^n-1$, of the linearized sequence $(k$~is {\rm not} the numerical address here$)$.

  Then $[$we conjecture$]$

  $(i)$ $p_{k+1}q_k - p_kq_{k+1}=2^e$ for some $e\in\nn_0$, $1\leq k\leq  2^n-2$.

  $(ii)$ The exponent $e$ is zero, the value thus 1, except for the following cases$:$

  \bt{llr|llr}
  Parent&Child&Value&Parent&Child&Value\\
  \cline{1-6}
  ${\tt B}_*$ or ${\tt C}_*$&${\tt C}_e$&$2^e$&
  $\overline{\tt B}_*$ or $\overline{\tt C}_*$&$\overline{\tt C}_e$&$2^e$\\
  ${\tt B}_{e-1}$&${\tt B}_{e}$&$2^e$&      
  $\overline{\tt B}_{e-1}$&$\overline{\tt B}_{e}$&$2^e$\\
  ${\tt A}$&${\tt B}_1$&$2$\ \ &
  $\overline{\tt A}$&$\overline{\tt B}_1$&$2$\ \ \\
  \et

\noindent  where the states are taken from Appendix~3.
\end{theorem}

\begin{proof}
  By numerical verification up to level $n=30$, {\it i.e.} for $1\leq k \leq 2^{30}-1$.
\end{proof}

\begin{theorem} Values as weighted Mediants between Neighbour Nodes 
 
  Let ${p_k}/{q_k}$ as in Theorem~\ref{pqpq} and 
  $$\Delta^+ = p_{k+1}q_k - p_kq_{k+1},\quad \Delta^- = p_kq_{k-1} - p_{k-1}q_k.$$
  Let $g = \operatorname{gcd}(\Delta^+,\Delta^-), \Delta_+ = \Delta^+/g, \Delta_- = \Delta^-/g$.
  Then
  $$\frac{p_k}{q_k} = \frac{\Delta_+ \cdot p_{k-1} + \Delta_-\cdot p_{k+1}}{\Delta_+ \cdot q_{k-1} + \Delta_-\cdot q_{k+1}}.$$
\end{theorem}

\begin{proof}
  \begin{eqnarray*}
    \frac{p_k}{q_k} &=& \frac{p_{k-1}\Delta^+ + p_{k+1}\Delta^-}{q_{k-1}\Delta^+ + q_{k+1}\Delta^-}\\
    &\Leftrightarrow& \left[q_{k-1}(p_{k+1}q_k-p_kq_{k+1}) + q_{k+1}(p_{k}q_{k-1}-p_{k-1}q_{k})\right]p_k\\
 &=&\left[p_{k-1}(p_{k+1}q_k-p_kq_{k+1}) + p_{k+1}(p_{k}q_{k-1}-p_{k-1}q_{k})\right]q_k\\  
   &\Leftrightarrow& p_kq_k \left[q_{k-1}p_{k+1}-q_{k+1}p_{k-1}\right]\\  
&=& p_kq_k\left[-p_{k-1}q_{k+1} + p_{k+1}q_{k-1}\right]
\end{eqnarray*}
  \end{proof}

\begin{note}
  The equivalent result for the Stern-Brocot tree is $e=0, \Delta_+=\Delta_- = 1$ for all nodes.
\end{note}

\begin{definition}

  For a given binary tree T with labels $p/q\in \qq$ at address $v(p/q)\in A^*$:
  
  $(i)$ Let $\lambda_T(q) =\ds\max_{\scriptsize \ba{c}1\leq p < q\\ (p,q)=1\ea} |v(p/q)|/\log_2(q)$ for the last level, such that all irreducible fractions with denominator $q$ are present in levels 1 to $\lambda_T(q)\cdot \log_2(q)$.

  $(ii)$ Let $\Lambda_T = \lim\sup _{q\in\nn} (\lambda_T(q))$.

\end{definition}

\begin{theorem} $($Almost$)$ Optimality of the V$_{10}$ Tree

  $(i)$ $\Lambda_T \geq 2$ for any tree T.

  $(ii)$ $\Lambda_{V_{10}} \geq 2.4007$.

  $(iii)$ $\Lambda_{V_{10}} \leq 3.44$.

  $(iv)$ For the Stern-Brocot tree, $\lambda_{SB}(q) \geq \frac{q}{\log_2(q)}$, and thus $\Lambda_{SB} = +\infty$.
  
\end{theorem}

\begin{proof}
$(i)$ There are $\phi(q)$ reduced fractions $p/q, 1\leq p < q$ in $\qq_1$.
Asymptotically, we have  $\sum_{k=1}^q \phi(k) \approx \frac{1}{2\zeta(2)} \cdot q^2 = \frac{3}{\pi^2} \cdot q^2$ values with denominator $\leq q$.
Therefore, for any  binary  tree we can at best expect to see all quotients with denominators $\leq  q$ in the first $2\log_2(q)$ levels, and thus $\Lambda_T \geq 2$ for any tree $T$.

$(ii)$ The irrational number $r=\sqrt{5}-2 = 0.236\dots$ with CFE $[4,4,4,\dots]$  has convergents $A_k/B_k$ with asymptotical growth of the denominator $B_k =  \Theta(\varphi_4^k) = (4.236\dots)^k$, and an encoding of $5k$ bits for the first $k$ copies of $b_i=4$.

Hence, the denominator $q=B_k$ appears (approximately, asymptotically) on level $5k$ in the V tree, where $(\log_2(q)\cdot\alpha\approx)\   \log_2(4.236^k)\cdot \alpha = 5k  \Leftrightarrow \alpha = 5/\log_2(2+\sqrt{5}) = 2.4007\dots$.

$(iii)$ We need at most $\log(q)/\log(\varphi_1)$ PDs at all, even if they all should be equal to 1.

Also, $q \leq \prod_i b_i$.
We advance in the product by a factor of 2, and 3 coding bits, or faster for other factors:
$\log(b) / l_{I,II}(b)$ is minimal for $b=2$ (except $b=1$, of course).
Hence, we get to the full product $q$ with at most  $3\cdot \log_2(q)$ coding bits for the PDs greater than 1, and at most
$1\cdot (\log_{\varphi_1}(q) - \log_2(q))$ bits for additional PDs with value 1 (which do not improve the product, but add to the coding length).
Hence, $\log_2(q) \times (2 + 1/\log_2(\varphi_1)) = 3.44\log_2(q)$ is the last level, where a denominator $q$ might appear.

$(iv)$ For the Stern-Brocot tree, $\lambda_{SB}(q) \geq \frac{q}{\log_2(q)}$, since  $1/q$ is on level $q$. $\Lambda_{SB}$ follows.
\end{proof}

\begin{note}
  $(i)$ Numerical evidence suggests $\Lambda_V\approx 2.5$.

  $(ii)$ Moving the lower bound  for $\Lambda$ below $2.35931 = 9/\log_2(\varphi_{14}) = 3/\log_2(\varphi_2)$  (the coincidence stems from $\varphi_{14} = \varphi_2^3$) is impossible with integral wordlengths, since then already $\sum_{b=1}^{64} 2^{-l(b)} > 1$.
  Hence, our encoding is basically optimal, besides being very regular.
\end{note}

\begin{conjecture}
\label{arc}
  Let $f$ be any continuous and monotonically increasing function $f\colon [0,1]\to[0,1]$ with $f(0)=0$ and $f(1)=1$.
  Then the graph $\{(x,f(x))\ |\ x\in[0,1]\}$ has Hausdorff dimension 1 and arc length between $\sqrt{2}$ and 2.
\end{conjecture}

\begin{proof} (idea)\quad 
We cover the graph by squares of side length $2^{-k}$, for $k\to\infty$, to show the upper bound  and the Hausdorff dimension.
The lower bound follows from the triangle inequality.
\end{proof}

Now, we will state some conjectures about the  graph of $?_V$ and $?_V^{-1}$.

\begin{conjecture}
 Assuming that the function $?_V(x)$ is continuous and monotonically increasing from 
$?_V(0)=0$ to $?_V(1)=1$, we conjecture that it has Hausdorff dimension 1, and in particular is not fractal.
\end{conjecture}

\begin{proposition}
\label{int}
  The area between the functions and the diagonal on $[0,1]$ satisfies
  $$0.030734101 < \int_{0}^1 (?_V(x) -x) \operatorname{dx} = -\int_{0}^1 (?^{-1}_V(v) -v) \operatorname{dv} < 0.030734102$$ 
\end{proposition}

\begin{proof} \quad 
The number 0.030734101$\dots$ results by taking the ``Riemann sum'' for the $2^k+1$ arguments $a/2^k, 0\leq a \leq  2^k$, for $k=1,\dots,30$.
The values settle.
\end{proof}

\begin{proposition}
  The arc length between $(0,0)$ and $(1,1)$ is greater than $1.554$.
\end{proposition}

\begin{proof}
  The lower bound for the arc length results by taking a polygonal chain  through the points $(x,?_V  (x))$ for the $2^k+1$ values $a/2^k, 0\leq a\leq 2^k$,  for $k=1,\dots,30$,  by summing up the length of the polygonal chain ($2^{30}$ diagonals).
\end{proof}

\begin{note}
  Since any finite number of points is compatible with the upper bound 2, by assuming that the curve also goes through the points $(x_{k+1}-\varepsilon,f(x_k)),\forall k$ for an arbitrarily small $\varepsilon > 0$, we can not improve that upper bound 2  for the arc length (from Conjecture~\ref{arc}) in this way.
\end{note}

\begin{conjecture} For all $x \in \qq_1$ (and by continuity in $\rr_1$), we have
$$\frac{8}{9} x \leq ?^{-1}(x) \leq x \leq ?(x)\leq \frac{9}{8} x$$
\end{conjecture}

Proof idea: $?_V(x)-x$ has minima for $?_V(2^{-k}) = 2^{-k}$ and 
maxima for $?_V(\frac{2}{3}\cdot 2^{-k}) = \frac{3}{4}\cdot 2^{-k}$,
and thus  $?_V^{-1}(x)-x$ has
maxima for $?_V^{-1}(2^{-k}) = 2^{-k}$ and 
minima for $?_V^{-1}(\frac{3}{4}\cdot 2^{-k}) = \frac{2}{3}\cdot 2^{-k}$.
This is, however, only verified numerically on points $a/2^k, k\leq 30$ from $\dd_1$.

\begin{note}
Self-similarity of the graph of the function $?(x)$ (see  Figure~\ref{fig:Graph})

While not fractal, the graph nevertheless exhibits a clear self-similarity:
$$?_V\left(\frac{x}{2}\right)\approx \frac{1}{2} ?_V(x), \forall x\in [0,1]$$

The dashed lines are the identity $y=x$ and $y=0$, respectively.
The dotted lines touch the local maxima $y = 9/8x$ and $y=1/8x \cdot 12$, respectively.
\end{note}

\begin{conjecture} Parabola Conjecture
  
  Apparently, in particular visible for $k=4$ in red in  Figure~\ref{fig:Graph}, the function graph is upper-bounded by curves through $(2^{-k},0)$ and $(2/3\cdot 2^{-k}, 3/4\cdot 2^{-k})$, which actually seem to be  parabolas for the inverse function $?_V^{-1}$.
  We thus conjecture:
  
  For $y\in\rr_1, \exists k\in\nn_0$ with $ \frac{1}{2}\cdot 2^{-k} \leq y \leq {1}\cdot 2^{-k}$.
  Using this $k$, we conjecture
  $$\overline ?_V^{-1}(y) \geq \frac{4}{3}2^k (y-\frac{3}{4}\cdot 2^{-k})^2 + y -\frac{1}{12}\cdot 2^{-k}$$
  which is met with equality (only) in  the three points $\frac{1}{2}\cdot 2^{-k},\ \frac{3}{4}\cdot 2^{-k},\ {1}\cdot 2^{-k}$.
  
  From $\int_{-1/4}^{1/4} \left(\frac{4}{3}y^2-\frac{1}{12}\right)\operatorname{dy} = -\frac{1}{36}$ and with $1+1/4+1/16+\dots = 4/3$, we have a combined  area of $\frac{1}{27}\approx 0.037$,  to be compared with the result 0.0307 from Proposition~\ref{int}.
\end{conjecture}

\newpage
\vspace*{-10mm}
\begin{figure}[!h]
\vspace*{-10mm}
  \begin{center}

  
  $k:$ 0 = black, {\color{green}1 = green}, {\color{yellow}2 = yellow}, {\color{blue}3 = blue}, {\color{red}4 = red}

\end{center}
\caption{\label{fig:Graph} Function graph and distance from diagonal.}
\vspace*{-20mm}
\end{figure}
\vspace*{-20mm}
\newpage

\begin{conjecture} The derivative  $\overline ?'_V(x)$ does not exist on $\qq_1$

  For $p/q\in\qq_1$, $p/q = [b_1,\dots,b_l] = [b_1,\dots,b_l-1,1], b_l\geq 2$, we have -- if defined at all (see example below):
  $$  \lim_{x\to \left(\frac{p}{q}\right)^-}  \overline ?'_V(x)= q^2\cdot 2^{-\alpha_L-\sum_{i=1}^l l_{I,II}(b_i)}$$
  and
  $$  \lim_{x\to \left(\frac{p}{q}\right)^+}  \overline ?'_V(x)= q^2\cdot 2^{-\alpha_R-\sum_{i=1}^l l_{I,II}(b_i)},$$
  where $\alpha_L\neq \alpha_R$ depend on $l$ and $b_l$:

  \bt{cc|cc|l}
  $l$&$b_l$&$\alpha_L$&$\alpha_R$\\
  \cline{1-5}
  odd& $\neq 2^k$ &  0 &+1   &additionally $b_{l+1}=1$ with $l_{I,II}(1)=1$ \\
  odd& $2^k$      &  0 &$-$1 & $l_{I,II}(2^k-1) = l_{I,II}(2^k)-2$, plus $b_{l+1}=1$ \\
  even& $\neq 2^k$& +1 &0    & as above, with sides reversed \\
  even& $2^k$     &$-$1&0    &\\ 
  \et

  Therefore, the derivative of  $\overline ?_V(x)$  does {\it not} exist at least in rational points (for Minkowski's $?(x)$, we have  $?'(x)=0$ for rational $x$, see \cite{Dush}). 
\end{conjecture}

\begin{example}
  $x_n= \overline ?_V^{-1}(?_V(38/51)\pm 2^{-n}), (38/51) = [1,2,1,12] = [1,2,1,11,1]$.
  \begin{eqnarray*}
    \lim_{x\to\left(\frac{38}{51}\right)^-}\overline ?'_V(x) &\stackrel{?}{=}&
    \lim_{n\to \infty}\frac{\overline ?_V\left(\frac{38\cdot 2^n+35}{51\cdot 2^n+47}\right)-\overline ?_V\left(\frac{38}{51}\right)}{\frac{38\cdot 2^n+35}{51\cdot 2^n+47} - \frac{38}{51}}\\
    &=&\lim_{n\to \infty}\frac{.1|100|1|1110011|1|1^n00^n|0^\omega - .1|100|1|1110100|0^\omega}
        {\frac{38\cdot 2^n+35}{51\cdot 2^n+47} -\frac{38}{51}}\\
        &=& \lim_{n\to \infty}
        \frac{\frac{829}{1024} -  2^{-n-\framebox{\scriptsize{13}}} -\frac{829}{1024}}
             {\frac{51\cdot 38\cdot 2^n+51\cdot 35-51\cdot 38\cdot 2^n -38\cdot 47}{51\cdot 51\cdot 2^n+ 47\cdot 51}}\\
             &=& \lim_{n\to \infty} %
             \frac{-(51^2\cdot 2^n+47\cdot 51)\cdot 2^{-n- 13}}{51\cdot 38\cdot 2^n -38\cdot 51\cdot 2^n+(51\cdot 3-38\cdot 47)}\\%
             &=&\lim_{n\to \infty} (-51^2\cdot 2^{-13} - 2^{-n}\cdot 47\cdot 51\cdot 2^{-13}) / (-1)
             = \frac{51^2}{2^{\framebox{\scriptsize{13}}}}
  \end{eqnarray*}
  \begin{eqnarray*}
    \lim_{x\to\left(\frac{38}{51}\right)^+} \overline ?'_V(x)&\stackrel{?}{=}&
    \lim_{n\to \infty}\frac{\overline ?_V\left(\frac{38\cdot 2^n+3}{51\cdot 2^n+4}\right)-\overline ?_V\left(\frac{38}{51}\right)}{\frac{38\cdot 2^n+3}{51\cdot 2^n+4} - \frac{38}{51}}\\
    &=&\lim_{n\to \infty}\frac{.1|100|1|1110100|0^n11^n|1^\omega - .1|100|1|1110100|0^\omega}
        {\frac{38\cdot 2^n+3}{51\cdot 2^n+4} -\frac{38}{51}}\\
        &=& \lim_{n\to \infty}
        \frac{\frac{829}{1024} +2^{-n-\framebox{\scriptsize{12}}} -\frac{829}{1024}}{\frac{51\cdot 38\cdot 2^n+51\cdot 3-51\cdot 38\cdot 2^n -38\cdot 4}{51\cdot 51\cdot 2^n+ 4\cdot 51}}\\
        &=& \lim_{n\to \infty}
        \frac{(51^2\cdot 2^n+4\cdot 51)\cdot 2^{-n-{\bf 12}}}{51\cdot 38\cdot 2^n -38\cdot 51\cdot 2^n+(51\cdot 3-38\cdot 4)}\\
        &=&\lim_{n\to \infty} (51^2\cdot 2^{-12} + 2^{-n}\cdot 4\cdot 51\cdot 2^{-12})/1
        = \frac{51^2}{2^{\framebox{\scriptsize{12}}}}
  \end{eqnarray*}
\end{example}

\begin{note}
  From {\large{$\frac{q^2}{2^{|\rm Code|}} \approx\frac{2^{3.42371\cdot l}}{2^{3.507\cdot l}}$}} $= 2^{-0.084\cdot l}$, by L\'evy,
  for every 12 PDs we should need some $3.42371 \cdot 12 \approx 41$ bits, but actually  we need one more, namely  $3.507\cdot 12\approx 42$ bits.
  This one more bit every 41 bits is + 2.4\% (compare with Proposition~\ref{entro},  $H_{C_{I,II}}/H_{L\acute{e}vy} = 1.024\dots)$.
\end{note}

\begin{note}
  The longer the CFE becomes, the flatter the (one-sided) derivatives at $p/q$.
\end{note}

\subsection *{Acknowledgement }

I wish to thank  my wife and active proofreader (valuable comments, debugging rate higher than mine ;-)
M\'onica del Pilar, with all my heart,
\hfill {\tiny{\# $\to$ "!`"}}

\vfill

Valdivia (Chile), Asunci\'on de la Virgen, A.D.~MMXX

\vfill

\vfill

\newpage
\section*{Appendices}

{\bf 1. Continued Fraction Expansion}

Let $r\in\rr$.
Let $\lfloor r \rfloor\in\zz$ be the largest integer smaller than or equal to $r$, {\it e.g.}
$\lfloor 3.14\rfloor = 3, \lfloor -3.14\rfloor = -4$, and let  $\{r\} = r - \lfloor r\rfloor  \in[0,1)$ be the fractional part.
  {\it E.g.}  $\{ 3.14\} = 0.14, \{ -3.14\} = 0.86.$
  
  The continued fraction expansion of $r=: r_0$ is defined by its successive partial denominators $b_i$ as
  $b_0 := \lfloor r_0 \rfloor, r_i := \frac{1}{\{r_{i-1}\}} =  \frac{1}{r_{i-1}-\lfloor r_{i-1}\rfloor},$
  $b_i := \lfloor r_i\rfloor\in\nn$, for $i\in\nn$.
  The continued fraction for $r$ is then 
  $$r = b_0+\frac{1}{b_1+\frac{1}{b_2+\frac{1}{\dots}}}
  = b_0+\frac{\hfill 1 \hfill|}{|\hfill b_1 \hfill}+\frac{\hfill 1 \hfill|}{|\hfill b_2 \hfill}+\cdots =: [b_0;b_1,b_2,\dots]$$
  and the convergents $A_i/B_i$ to $r$ are obtained by  Perron's schema \cite[S.~24]{Perr} (Table~\ref{fig:CFE-Perron}).
The  initial values are $B_{-2} = A_{-1} = 1$, $A_{-2} = B_{-1} = 0$ and then $A_i := b_i\cdot A_{i-1}+A_{i-2}$, $B_i := b_i\cdot B_{i-1}+B_{i-2}$.
In particular $A_0 = b_0, B_0 = 1,
A_1 = b_1 A_0 + A_{-1} = b_1b_0+1,
B_1 = b_1 B_0 + B_{-1} = b_1$.
We focus on the case $r\in(0,1) =\rr_1\subset \rr$, thus $b_0=0$ ({\it e.g.}~for $r=\pi-3$ see the second part of Table~\ref{fig:CFE-Perron}).

\begin{table}[!h]
  \begin{center}
    {\small{
        $\ba{c|c|c|c|c|c|c|c|c|c}
        i  &-2&-1&  0&  1&  2&  3&  4&\dots\\
        b_i& -& -&b_0&b_1&b_2&b_3&b_4&\dots\\
        \cline{1-9}                    
        A_i& 0& 1&A_0&A_1&A_2&A_3&A_4&\dots\\
        B_i& 1& 0&B_0&B_1&B_2&B_3&B_4&\dots\\
        \cline{1-9}
        b_i& & & &7&15&1&292&$\dots$\\
        \cline{1-9}
        A_i&0&1&0&1& 15& 16&4786&$\dots$\\
        B_i&1&0&1&7&106&113&33102&$\dots$\\
        \cline{1-9}
        \ea$
    }}                             
  \end{center}                     
  \caption{\label{fig:CFE-Perron} CFE Schema according to Perron.}
\end{table}

\noindent Convergence:\quad 
For $r\in\rr^+$, we have\\
$0 = \frac{A_{-2}}{B_{-2}} \leq  \frac{A_{0}}{B_{0}} < \frac{A_{2}}{B_{2}} < \frac{A_{4}}{B_{4}} < \cdots < r < \cdots < \frac{A_{5}}{B_{5}} < \frac{A_{3}}{B_{3}} < \frac{A_{1}}{B_{1}} < \frac{A_{-1}}{B_{-1}} = \infty$\linebreak
and furthermore $\left|r - \frac{A_{k}}{B_{k}}\right| < \frac{1}{B_kB_{k+1}}$, \cite[Satz 2.10]{Perr}.
\\\\
Ambiguity:\quad 
$[b_1,\dots, b_l] = [b_1,\dots,b_l-1,1]$ and  $[b_1,\dots, b_l, \aleph_0] = [b_1,\dots,b_l-1,1, \aleph_0]$

One can  resolve this ambiguity in 4 ways:\\
\noindent $(i)$ Let the last PD be always greater than 1, or\\
$(ii)$ always equal to 1, or\\
$(iii)$ have an even, or \\
$(iv)$ an odd number of PDs\\
(the final $\aleph_0$  with $ 1/\aleph_0 := 0$ in any case does not alter the value).
\\\\
\noindent We shall use convention $(iii)$: The encoding then terminates in $0^\omega$ from $\Ci(\aleph_0) = 0^\omega$.

\newpage
{\bf 2. Equivalence between $\nn, A^*$, and $\dd_1$}

We identify the word $v\in A^*$ with the number $n=(1v)_2\in\nn_0$ in binary representation, and the dyadic fraction $(v|1)_2/2^{|v|+1}\in\dd_1$.

In particular:\\
$v=\varepsilon \equiv n=1 \equiv d = 1/2,\\
v=0 \equiv n=2 \equiv d = 1/4,\\
v=1 \equiv n=3 \equiv d = 3/4$.

Example:\\
$v = 10010$ (value 18) $\equiv\\
n = 1|10010_2 = 50 = 18+2^5\equiv\\
d = p/2k=(10010|1)/2^6 = 37/64$ where $37 = 18\cdot 2+1$.
\\\\
We define bijective mappings between the 3 sets  $\nn, A^*$, and $\dd_1$ as follows,\\
where $\iota_{XY}^{-1} = \iota_{YX}$ for $X,Y \in \{N,A,D\}$ and $l :=\lfloor \log_2(n)\rfloor$:

$\ba{cclll}
\\
\iota_{NA}\colon& \nn\to A^*,     & \iota_{NA}(n)     &=& n-2^l\mbox{\rm\  in\ binary}\\
\\
\iota_{AN}\colon& A^*\to \nn,     & \iota_{AN}(v)     &=& (v)_2 + n+2^{|v|+1}\\
\\
\iota_{ND}\colon& \nn\to \dd_1,   & \iota_{ND}(n)     &=& ((n-2^l)\cdot 2+1)/2^l\\
\\
\iota_{DN}\colon& \dd_1\to \nn,   & \iota_{DN}(p/2^k) &=& (p-1)/2 +2^{k-1}\\
\\
\iota_{AD}\colon& A^*\to\dd_1,    & \iota_{AD}(v)     &=& (v|1)_2/2^{|v|+1}\\
\\
\iota_{DA}\colon& \dd_1\to A^*,   & \iota_{DA}(p/2^k) &=& 0^{k-|p|-2}|((p-1)/2)\mbox{\rm\  in\ binary}\\
\ea$

\vfill

{\bf 3. Finite State Machine} …
\label{FSM}

Let $Q=\{A, B_k, C_k,\overline A, \overline B_k, \overline C_k, k\in\nn\}$ be the state set for an FSM with nextstate function
$Q^+\colon Q\times A\to Q$ given by:
\\\\
$\ba{l|l|l|l}
q&Q^+(q,0)&Q^+(q,1)\\
\cline{1-3}
A&B_1  & \overline A &\mbox{\rm \ Start\ for\ }  \Ci    \\
B_k&B_{k+1}  & C_k&\mbox{\rm \ Increase\ PDs as\ } 1,2,4,8,16,\dots\\
C_k,k\geq 2& C_{k-1} & C_{k-1}&\mbox{\rm \ Adjust\ PDs\ by\ }\pm ...16,8,4,2    \\
C_1&\overline A  & \overline A&\mbox{\rm \  Adjust\ PD\ by\  }\pm 1,\mbox{\rm \  switch\ to\ }\Cj\\
\overline A&A  & \overline B_1&\mbox{\rm \ Start\ for\ }  \Cj     \\
\overline B_k& \overline C_k & \overline B_{k+1} &\mbox{\rm \ Increase\ PDs as\ } 1,2,4,8,16,\dots\\
\overline C_k, k\geq 2& \overline C_{k-1} &  \overline C_{k-1}  &\mbox{\rm \ Adjust\ PDs\ by\ }\pm ...16,8,4,2    \\
\overline C_1&A  & A  &\mbox{\rm \  Adjust\ PD\ by\  }\pm 1,\mbox{\rm \  switch\ to\ }\Ci\\
\ea$

\vfill

\vfill

\newpage
{\bf 4. Trees and Addresses}

The numerical address  $n\in \nn$ and the symbolic address $v\in A^*$ are related by  $n = (1v)_2$ in binary, see Appendix~2.
{\it E.g.} on the last line we see  $n = $ {\bf 23} and  $v=0111$, with $23 = (1|0111)_2$
Note that the left child node has $n_L= 2n$ and $v_L=v0$,
the right one $n_R = 2n+1$ and  $v_R=v1$.
The dyadic fraction is $a/2^k\in \dd$, $a$ odd, with $a/2^k = (v1)_2/2^{|v|+1}$, and it comes from the van der Corput sequence in base 2 (see \cite[p.~127]{Kuip}), which is just $A^*$, the words written from right to left: $\varepsilon$,0,1,00,10,01,11,000,100,010,110,001,101,011,111,0001,\dots .

Here, the  dyadic fraction is  $(0111|1)_2/2^{|0111|+1} = 15/32$.
The three entries of the upper part coincide according to Appendix~2.

The bottom part consists of the two values from the Stern-Brocot tree and from the V$_{10}$  tree.
The entry here is 4/9 for both trees.
Using $\iota_{DA}$ from Appendix~2, we can say that
$?^{-1} \circ \iota_{DA}$  maps the van der Corput tree to the Stern-Brocot tree, and 
$?_V^{-1} \circ \iota_{DA}$  maps the van der Corput tree to the V$_{10}$  tree, entry by entry.

\begin{figure}[!h]
  \begin{center}
    {\scriptsize{
\hspace*{-2 cm}
        \begin{forest}
          [{\framebox{\hq$\ba{c}\mbox{\bf 1}\\ \varepsilon\\
                {1}/{2}\\
                \\
                1/2\\
                1/2\\
                \ea$\hq}}
            [{\framebox{\hq$\ba{c}\mbox{\bf 2}\\0 \\
                  {1}/{4}\\
                  \\
                  1/3\\
                  1/4\\
                  \ea$\hq}}
              [{\framebox{\hq$\ba{c}\mbox{\bf 4}\\00 \\
                    {1}/{8}\\
                    \\
                    1/4\\
                    1/8\\
                    \ea$\hq}}
                [{\framebox{\hq$\ba{c}\mbox{\bf 8}\\000 \\
                      {1}/{16}\\
                      \\
                      1/5\\
                      1/16\\
                      \ea$\hq}}
                  [{\hq\framebox{\hq$\ba{c}\mbox{\bf 16}\\0000 \\
                        {1}/{32}\\
                        \\
                        1/6\\
                        1/32\\
                        \ea$\hq}\hq}]
                  [{\hq\framebox{\hq$\ba{c}\mbox{\bf 17}\\0001 \\
                        {3}/{32}\\
                        \\
                        2/9\\
                        1/12\\
                        \ea$\hq}\hq}]
                ]
                [{\framebox{\hq$\ba{c}\mbox{\bf 9}\\001 \\
                      {3}/{16}\\
                      \\
                      2/7\\
                      1/6\\
                      \ea$\hq}}
                  [{\hq\framebox{\hq$\ba{c}\mbox{\bf 18}\\0010 \\
                        {5}/{32}\\
                        \\
                        3/11\\
                        1/7\\
                        \ea$\hq}\hq}]
                  [{\hq\framebox{\hq$\ba{c}\mbox{\bf 19}\\0011 \\
                        {7}/{32}\\
                        \\
                        3/10\\
                        1/5\\
                        \ea$\hq}\hq}]
                ]
              ]
              [{\framebox{\hq$\ba{c}\mbox{\bf 5}\\ 01 \\
                    {3}/{8}\\
                    \\
                    2/5\\
                    1/3\\
                    \ea$\hq}}
                [{\framebox{\hq$\ba{c}\mbox{\bf 10}\\010 \\
                      {5}/{16}\\
                      \\
                      3/8\\
                      2/7\\
                      \ea$\hq}}
                  [{\hq\framebox{\hq$\ba{c}\mbox{\bf 20}\\0100 \\
                        {9}/{32}\\
                        \\
                        4/11\\
                        3/11\\
                        \ea$\hq}\hq}]
                  [{\hq\framebox{\hq$\ba{c}\mbox{\bf 21}\\0101 \\
                        {11}/{32}\\
                        \\
                        5/13\\
                        4/13\\
                        \ea$\hq}\hq}]
                ]
                [{\framebox{\hq$\ba{c}\mbox{\bf 11}\\011 \\
                      {7}/{16}\\
                      \\
                      3/7\\
                      2/5\\
                      \ea$\hq}}
                  [{\hq\framebox{\hq$\ba{c}\mbox{\bf 22}\\0110\\
                        {13}/{32}\\
                        \\
                        5/12\\
                        3/8\\
                        \ea$\hq}\hq}]
                  [{\hq\framebox{\hq$\ba{c}\mbox{\bf 23}\\0111\\
                        {15}/{32}\\
                        \\
                        4/9\\
                        4/9\\
                        \ea$\hq}\hq}]
                ]
              ]
            ]
            [{\framebox{\hq$\ba{c}\mbox{\bf 3}\\1\\
                  {3}/{4}\\
                  \\
                  2/3\\
                  2/3\\
                  \ea$\hq}}
              [{\framebox{\hq$\ba{c}\mbox{\bf 6}\\10 \\
                    {5}/{8}\\
                    \\
                    3/5\\
                    3/5\\
                    \ea$\hq}}
                [{\framebox{\hq$\ba{c}\mbox{\bf 12}\\100 \\
                      {9}/{16}\\
                      \\
                      4/7\\
                      5/9\\
                      \ea$\hq}}
                  [{\hq\framebox{\hq$\ba{c}\mbox{\bf 24}\\1000 \\
                        {17}/{32}\\
                        \\
                        5/9\\
                        9/17\\
                        \ea$\hq}\hq}]
                  [{\hq\framebox{\hq$\ba{c}\mbox{\bf 25}\\1001 \\
                        {19}/{32}\\
                        \\
                        7/12\\
                        4/7\\
                        \ea$\hq}\hq}]
                ]
                [{\framebox{\hq$\ba{c}\mbox{\bf 13}\\101 \\
                      {11}/{16}\\
                      \\
                      5/8\\
                      5/8\\
                      \ea$\hq}}
                  [{\hq\framebox{\hq$\ba{c}\mbox{\bf 26}\\1010 \\
                        {21}/{32}\\
                        \\
                        8/13\\
                        8/13\\
                        \ea$\hq}\hq}]
                  [{\hq\framebox{\hq$\ba{c}\mbox{\bf 27}\\1011 \\
                        {23}/{32}\\
                        \\
                        7/11\\
                        9/14\\
                        \ea$\hq}\hq}]
                ]
              ]
              [{\framebox{\hq$\ba{c}\mbox{\bf 7}\\11\\
                    {7}/{8}\\
                    \\
                    3/4\\
                    4/5\\
                    \ea$\hq}}
                [{\framebox{\hq$\ba{c}\mbox{\bf 14}\\110 \\
                      {13}/{16}\\
                      \\
                      5/7\\
                      3/4\\
                      \ea$\hq}}
                  [{\hq\framebox{\hq$\ba{c}\mbox{\bf 28}\\1100\\
                        {25}/{32}\\
                        \\
                        7/10\\
                        5/7\\
                        \ea$\hq}\hq}]
                  [{\hq\framebox{\hq$\ba{c}\mbox{\bf 29}\\1101\\
                        {27}/{32}\\
                        \\
                        5/7\\
                        7/9\\
                        \ea$\hq}\hq}]
                ]
                [{\framebox{\hq$\ba{c}\mbox{\bf 15}\\ 111 \\
                      {15}/{16}\\
                      \\
                      4/5\\
                      8/9\\
                      \ea$\hq}}
                  [{\hq\framebox{\hq$\ba{c}\mbox{\bf 30}\\1110 \\
                        {29}/{32}\\
                        \\
                        7/9\\
                        6/7\\
                        \ea$\hq}\hq}]
                  [{\hq\framebox{\hq$\ba{c}\mbox{\bf 31}\\1111\\
                        {31}/{32}\\
                        \\
                        5/6\\
                        16/17\\
                        \ea$\hq}\hq}]
                ]
              ]
            ]
          ]
        \end{forest}
    }}                             
  \end{center}                     
  \caption{\label{fig:VDCTreeo} Addresses and trees: van der Corput, Stern-Brocot, and V$_{10}$  tree.}
\end{figure}
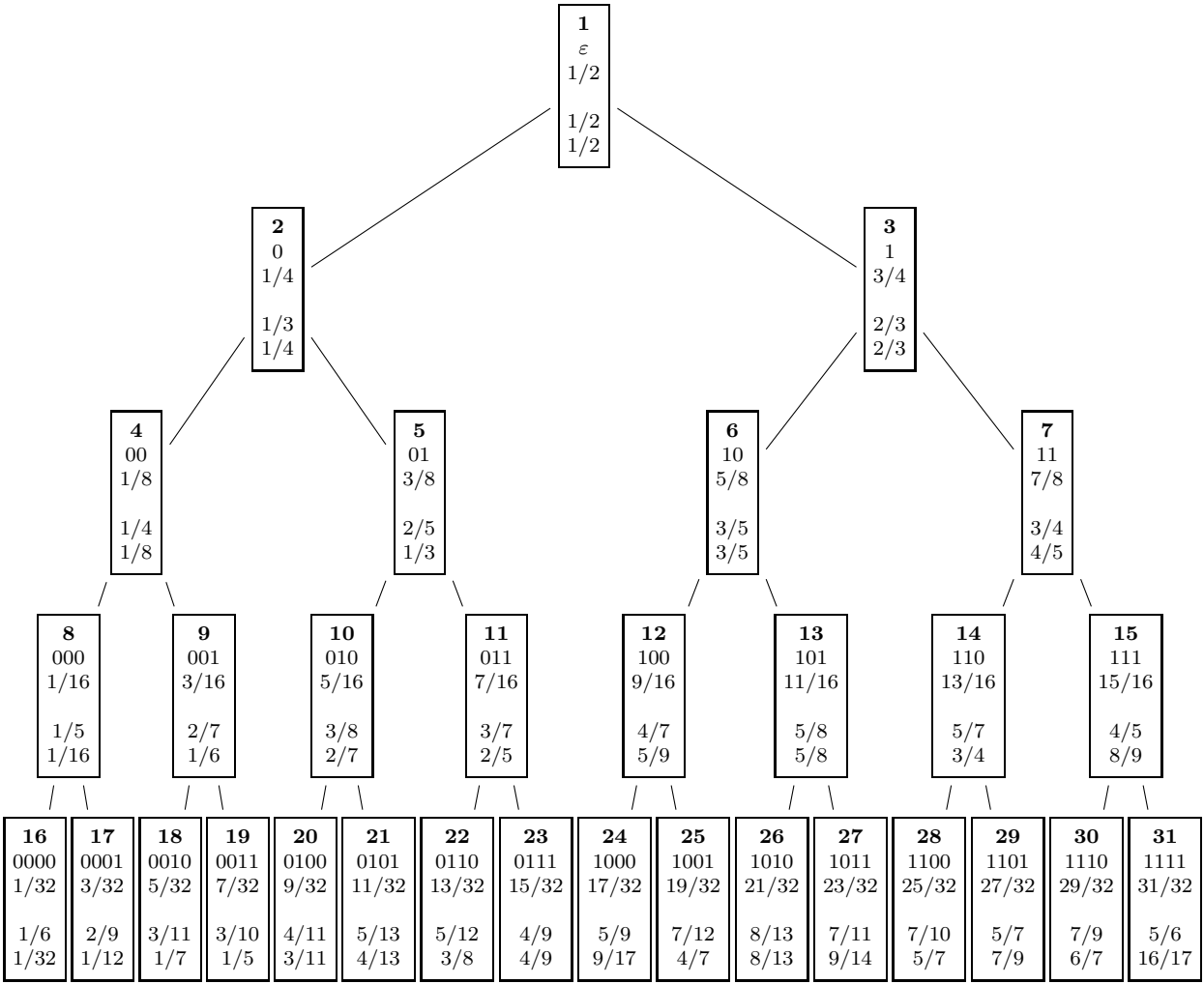

\newpage
{\bf 5. V$_{10}$ Values and their PDs and Encodings, for $|v|\leq 5$}

\begin{table}[!h]
  \begin{center}
    {\small{
        $\ba{llcl|llc}                                              
        v.10^\infty&\mbox{PDs}&A/B&r=(A/B)_2                    &   v.10^\infty&\mbox{PDs}&A/B\\                            
        \cline{1-7}       
        \varepsilon.1|0|&1|1|&1/2&.1(0)                                    & 00000.100000|0|\hs& 63|1|     & 1/64\\             
        \cline{1-4}       
        0.10|0|&    3|1|&1/4& .01(0)                            & 00001.1000|0|  & 23|1|     & 1/24\\               
        |1.00|&    1|2|&2/3& .(10)                              & 00010.10|0|    & 13|1|     & 1/14\\                 
        \cline{1-4}       
        00.100|0|&  7|1|&1/8& .001(0)                           & 00011.10|0|    & 9|1|          & 1/10\\                 
        01.1|0|&    2|1|&1/3& .(01)                             & 00100|.100|    & 7|2|          & 2/15\\                
        1|0|.1|0|&  1|1|1|1|&3/5& .(1001)                       & 00101|.100|    & 6|2|          & 2/13\\            
        1|1.1000|&  1|4|&4/5& .(1100)                           & 00110|.100|    & 5|2|          & 2/11\\            
        \cline{1-4}       
        000.1000|0|& 15|1|&1/16& .0001(0)                       & 00111|.100|    & 4|2|          & 2/9\\            
        001.10|0|&   5|1|&1/6&.00(10)                           & 010|0|0.10|0|  & 3|1|2|1|      & 4/15\\            
        010|.100|&   3|2|&2/7& .(010)                           & 010|0|1|.100|  & 3|1|1|2|      & 7/25\\            
        011|.100|&   2|2|&2/5&.(0110)                           & 010|10.1|      & 3|3|          & 3/10\\                     
        1|0|0.10|0|& 1|1|3|1|&5/9& .(100011)                    & 010|11.10000|  & 3|8|          & 8/25\\            
        1|0|1|.100|& 1|1|1|2|&5/8& .101(0)                      & 011|0|0.10|0|  & 2|1|3|1|      & 4/11\\        
        1|10.1|&    1|3|&     3/4& .11(0)                       & 011|0|1.|100|  & 2|1|1|2|      & 5/13\\        
        1|11.10000|& 1|8|&8/9&.(111000)                         & 011|10.1|      & 2|3|          & 3/7\\        
        \cline{1-4}       
        0000.10000|0|\hs\hs& 31|1|&1/32& .00001(0)              & 011|11.10000|  & 2|8|          & 8/17\\         
        0001.100|0|& 11|1|&1/12& .00(01)                        & 1|0|000.1000|0|\hs& 1|1|8|1|   & 10/19\\        
        0010.1|0|&   6|1|&1/7&.(001)                            & 1|0|001.10|0|  & 1|1|5|1|      & 7/13\\        
        0011.1|0|&   4|1|&1/5&.(0011)                           & 1|0|010|.100|  & 1|1|3|2|      & 9/16\\       
        010|0.|1|0|& 3|1|1|1|&3/11&.(0100010111)                & 1|0|011|.100|  & 1|1|2|2|      & 7/12\\       
        010|1.1000|& 3|4|&4/13& .(010011101100)                 & 1|0|1|0|0.10|0|    & 1|1|1|1|3|1|\hs& 14/23\\ 
        011|0|.1|0|& 2|1|1|1|&3/8&.011(0)                       & 1|0|1|0|1.|100|    & 1|1|1|1|1|2|\hs& 13/21\\      
        011|1.1000|& 2|4|&4/9& .(011100)                        & 1|0|1|10.1|        & 1|1|1|3|     & 7/11\\ 
        1|0|00.100|0|& 1|1|7|1|&9/17&.(10000111)                & 1|0|1|11.100|      & 1|1|1|8|     & 17/26\\       
        1|0|01.1|0|& 1|1|2|1|&4/7& .(100)                       & 1|100|0.10|0|      & 1|2|3|1|     & 9/14\\       
        1|0|1|0|.1|0|& 1|1|1|1|1|1|\hs&8/13& .(100111011000)    & 1|100|1|.100|      & 1|2|1|2|     & 8/11\\       
        1|0|1|1.1000|& 1|1|1|4|&9/14& .1(010)                   & 1|101|0.10|0|      & 1|3|3|1|     & 13/17\\ 
        1|100|.1|0|&1|2|1|1|& 5/7& .(101)                       & 1|101|1.|100|      & 1|3|1|2|     & 11/14\\
        1|101|.1|0|&1|3|1|1|& 7/9& .(110001)                    & 1|1100.1|          & 1|5|         & 5/6\\  
        1|110.10|  &1|6|    &6/7 &.(110)                        & 1|1101.1|          & 1|7|         & 7/8\\  
        1|111.100000|&1|16|  &16/17&.(11110000)                 & 1|1110.100|        & 1|12|        & 12/13\\  
        &&&                                                     & 1|1111.10000| \hs  & 1|32|        & 32/33\\  
        \ea$
    }}                             
  \end{center}                     
  \caption{\label{fig:CFE-Apx} Binary CFE and approximations.}
\end{table}

\end{document}